\input amstex
\input epsf
\documentstyle{amsppt}
%\magnification 1200
%\hsize=6truein
%\vsize=8.5truein

\def\ls{\leqslant}
\def\gs{\geqslant}
\def\SS{\frak S}

\TagsOnRight

\topmatter
\title
Layered restrictions and Chebyshev polynomials
\endtitle
\author Toufik Mansour$^*$ and Alek Vainshtein$^\dag$ \endauthor
\affil $^*$ Department of Mathematics\\
$^\dag$ Department of Mathematics and Department of Computer Science\\ 
University of Haifa, Haifa, Israel 31905\\ 
{\tt tmansur\@study.haifa.ac.il},
{\tt alek\@mathcs.haifa.ac.il}
\endaffil

\abstract
A permutation is called layered
if it consists of the disjoint union of substrings (layers)
so that the entries decrease within each layer, and 
increase between the layers. We find the
generating function for the number of permutations on $n$ letters avoiding
$(1,2,3)$ and a layered permutation on $k$ letters. In the most interesting
case of two layers, the generating function depends only on $k$ and is
expressed via Chebyshev polynomials of the second kind.
\medskip
\noindent {\smc 2000 Mathematics Subject Classification}: 
Primary 05A05, 05A15; Secondary 30B70, 42C05
\endabstract

\rightheadtext{Layered restrictions and Chebyshev polynomials}
\leftheadtext{Toufik Mansour and Alek Vainshtein}
\endtopmatter

%\vfill\eject

\document
\heading 1. Introduction and Main Result\endheading

Let $\pi\in\SS_n$ and $\tau\in\SS_k$ be two permutations.
An {\it occurrence\/} of $\tau$ in $\pi$
is a subsequence $1\ls i_1<i_2<\dots<i_k\ls n$ such that $(\pi_{i_1},
\dots,\pi_{i_k})$ is order-isomorphic to $\tau$; in such a context $\tau$ is
usually called a {\it pattern\/}. We say that $\pi$ {\it avoids\/} $\tau$,
or is $\tau$-{\it avoiding\/}, if there is no occurrence of $\tau$ in $\pi$.
The set of all $\tau$-avoiding permutations in $\SS_n$ is denoted 
$\SS_n(\tau)$. 
For an arbitrary finite collection of patterns $T$, we say that $\pi$
avoids $T$ if $\pi$ avoids any $\tau\in T$; the corresponding subset of $\SS_n$
is denoted $\SS_n(T)$.

In several recent papers \cite{CW, MV, Kr}, the authors express generating
functions for the number of permutations avoiding certain collections $T$ in
terms of Chebyshev polynomials. More specifically, the main result of
\cite{CW} can be formulated as follows. 

\proclaim{Theorem 1.1 ({\rm \cite{CW, Theorem~3.1})}}
Let $T_1=\{(1,2,3), 
(k-1,k-2,\dots,1,k)\}$, $T_2=\{(2,1,3), (1,2,\dots,k)\}$, and
$T_3=\{(2,1,3), (k,1,2,\dots,k-1)\}$, then the sequences $|\SS_n(T_i)|$,
$i=1,2,3$, have the same generating function
$$
R_k(x)=\frac{U_{k-1}\left(\frac1{2\sqrt{x}}\right)}{\sqrt{x}
U_{k}\left(\frac1{2\sqrt{x}}\right)},
$$
where $U_r(\cos\theta)={\sin(r+1)\theta}/{\sin\theta}$ is the Chebyshev 
polynomial of the second kind.
\endproclaim

Paper \cite{MV} gives a different proof of Theorem~1.1 for the case of $T_2$,
and provides explicit expressions for the generating functions of permutations
avoiding $(2,1,3)$ and containing $(1,2,\dots,k)$ a fixed number of times,
or avoiding  $(1,2,\dots,k)$ and containing $(2,1,3)$, via the same Chebyshev
polynomials. Still another proof
of Theorem~1.1 (for all the three cases) based on the use of Dyck paths is
given in \cite{Kr}, together with the corresponding generalizations
(to permutations avoiding the shorter pattern and containing the longer one a 
fixed number of times). 

In this paper, we give a completely different generalization of Theorem~1.1
for the case of $T_1$. Following \cite{Bo}, we say that $\tau\in \SS_k$ is 
{\it $p$-layered\/} if it consists of the disjoint union of $p$ substrings 
(the {\it layers\/}) so that the entries decrease within each layer, and 
increase between the layers. For any $k_1,\dots,k_p\gs1$ we
denote by $w(k_1,\dots,k_p)$ the unique layered pattern of length
$k=k_1+\dots+k_p$ whose $i$th layer is of length $k_i$. For example, the
only $1$-layered pattern of length $k$ is $w(k)=(k,k-1,\dots,1)$, and the 
only $k$-layered pattern of length $k$ is $w(1,\dots,1)=(1,2,\dots,k)$.
Thus, the longer pattern in the collection $T_1$ above is just the $2$-layered
pattern $w(k-1,1)$. 

The main result of this paper is as follows.

\proclaim{Theorem 1.2}
Let $T=\{(1,2,3),\tau\}$, where $\tau\in \SS_k$, $k\gs p$,
is an arbitrary $p$-layered 
pattern, and let $F_T(x)=\sum_{n\gs0}|\SS_n(T)|x^n$. Then\/{\rm:}

{\rm (i)} if $p=1$ then $F_T(x)$ is a polynomial of degree $2k-2${\rm;}

{\rm (ii)} if $p=2$ then $F_T(x)=R_k(x)${\rm;}

{\rm (iii)} if $p\gs 3$ then $F_T(x)=(1-\sqrt{1-4x})/2x$.
\endproclaim

In fact, only statement (ii) is new, while (i) and (iii) are added 
for the sake of completeness.
%in order to complete the classification. 
Indeed, if $p=1$ then $\tau=w(k)$, 
and by \cite{ES}, there exists a permutation in $\SS_n$
avoiding both $(1,2,\dots,l)$ and $(m,m-1,\dots,1)$ for
$n= (m-1)(l-1)$, and there are no such permutations for any $n> (m-1)(l-1)$. 
On the other hand, if $p\gs3$ then $\tau$ evidently
contains $(1,2,3)$, hence $|\SS_n(T)|=|\SS_n(1,2,3)|=c_n$, the $n$th Catalan 
number (see e.g. \cite{Kn, SS}), and the generating function of Catalan numbers
is known to be $(1-\sqrt{1-4x})/2x$.

Our proof of Theorem~1.2 is based on finding a recursion for the numbers
in question by purely analytical means. In spite of the paradigm
formulated in \cite{Kr}, that any enumeration problem leading to
Chebyshev polynomials is related to Dyck paths, it would be tempting
to find a proof that exploits such a relation. 

The final version of this paper was written during the second author's
(A.V.) stay at Max--Planck--Institut f\"ur Mathematik in Bonn,
Germany. A.V.~wants to express his gratitude to MPIM for the support.

\heading 2. Proofs \endheading

Denote by $T^{d,k}$ the collection $\{(1,2,3),w(k-d,d)\}$, and by $f^{d,k}_n$
the number of permutations in $\SS_n(T^{d,k})$, $1\ls d\ls k-1$. 
The following
observation is obtained straightforwardly from the definitions and usual
symmetry operations (reversal and complement).

\proclaim{Lemma 2.1} $f^{d,k}_n=f^{k-d,k}_n$ for all $d$, $1\ls d\ls k-1$.
\endproclaim

In what 
follows we assume that $k\gs 3$ and $d$ are fixed, and
omit the indices $d$, $k$ whenever appropriate; for example, instead of
$f^{d,k}_n$ we write just $f_n$. Moreover, it follows from Lemma~2.1
that we may assume that 
$$
k-d\gs d.  \tag1
$$  

For any $n\gs 1$ and any $m$ such that $1\ls m\ls n$, we denote by 
$g_n(i_1,\dots,i_m)=g_{n}^{d,k}(i_1,\dots,i_m)$ the number of permutations
$\pi\in \SS_n(T^{d,k})$ such that $\pi_j=i_j$ for $j=1,\dots,m$;
the corresponding subset of $\SS_n(T^{d,k})$ is denoted by 
$G_n(i_1,\dots,i_m)$. 

The following properties of the numbers $g_n(i_1,\dots,i_m)$ can be deduced
easily from the definitions.

\proclaim{Lemma 2.2} {\rm (i)} Let $n\gs 2$ and $1\ls i\ls n$, then
$$
g_n(\dots,i,\dots,i,\dots)=0.
$$

{\rm (ii)} Let $n\gs 3$, $1\ls m\ls n-2$, and $n-2\gs i_1>i_2>\dots>i_m\gs1$,
then
$$
g_n(i_1,\dots,i_m,j)=0\quad\text{for $\quad i_m+1\ls j\ls n-1$}.
$$ 

{\rm (iii)} Let $n\gs 3$, $1\ls m\ls \min\{n, k-d\}-1$, and 
$n-1\gs i_1>i_2>\dots>i_m\gs1$, then
$$
g_n(i_1,\dots,i_m,n)=g_{n-1}(i_1,\dots,i_m).
$$

{\rm (iv)}  Let $n\gs 3$, $2\ls m\ls n$, $i_1\gs n-d+1$, and 
$n\gs i_1>i_2>\dots>i_m\gs1$, then
$$
g_n(i_1,\dots,i_m)=g_{n-1}(i_2,\dots,i_m).
$$

{\rm (v)}  Let $n\gs k$, $k-d\ls m\ls n-d$, and 
$n-d\gs i_1>i_2>\dots>i_m\gs1$, then
$$
g_n(i_1,\dots,i_m)=0.
$$
\endproclaim

\demo{Proof} Property (i) is evident. To prove (ii) is enough to observe
that if $\pi\in G_n(i_1,\dots,i_m,j)$ then the entries $i_m$, $j$, and $n$
give an occurrence of $(1,2,3)$ in $\pi$, a contradiction. 

To prove (iii),
denote by $\pi^*$ the permutation obtained from $\pi$ by deleting its largest
entry. Then $\pi\in G_n(i_1,\dots,i_m,n)$ if and only if 
$\pi^*\in G_{n-1}(i_1,\dots,i_m)$, since entry $n$ placed as in (iii) cannot
be used in an occurrence of $(1,2,3)$ or $w(k-d,d)$ in $\pi$. 

The proof of
(iv) goes along the same lines. We denote by $\pi^\circ$ the permutation
obtained from $\pi$ by deleting entry $i_1$ and subtracting $1$ from all
the entries greater than $i_1$; then $\pi\in G_n(i_1,\dots,i_m)$ if and 
only if $\pi^\circ\in G_{n-1}(i_2,\dots,i_m)$. Indeed, the only if part is
trivial, so assume that $\pi^\circ\in G_{n-1}(i_2,\dots,i_m)$, and there 
exists an occurrence of $(1,2,3)$ or $w(k-d,d)$ in $\pi$. Evidently, such
an occurrence makes use of $i_1$. If $(i_1,j_1,j_2)$ is an occurrence
of $(1,2,3)$, then $j_1$ and $j_2$ lie to the right of $i_m$, and since 
$i_1>i_m$, the triple $(i_m,j_1,j_2)$ is an occurrence of $(1,2,3)$ in
$\pi^\circ$, a contradiction. Otherwise, $i_1$ enters an occurrence of
$w(k-d,d)$, and hence $\pi$ contains at least $d$ entries strictly greater
than $i_1$, a contradiction. 

Finally, to prove (v), consider all the entries of $\pi$ greater than $i_1$.
There are at least $d$ such entries, and all of them lie to the right
of $i_m$. If a pair $(j_1,j_2)$ of these entries is an occurrence of $(1,2)$,
then $(i_m,j_1,j_2)$ is an occurrence of $(1,2,3)$, a contradiction. Hence,
all these entries provide an occurrence of $(p,p-1,\dots,1)$ for some 
$p\gs d$; therefore, any $d$ out of these entries, together with 
 any $k-d$ entries out of $i_1,\dots,i_m$, give an occurrence of $w(k-d,d)$.
\qed
\enddemo

Now we introduce several objects that play an important role in the proof of
the main result. For $n\gs d+1$ and $1\ls m\ls \min\{n, k\}-d$ we define
$$\gather
\Cal I(n,m)=\{(i_1,\dots,i_m)\: m-j\ls i_{j+1}\ls i_j-1, j=0,\dots,m-1\},
i_0=n-d+1,\\
\Cal J(n,m)=\{(i_1,\dots,i_m)\: m-j\ls i_{j+1}\ls i_j-1, j=0,\dots,m-1\},
i_0=n-d+2,\\
\Cal K(n,m)=\{(i_1,\dots,i_m)\: (i_1,\dots,i_{m-1})\in \Cal I(n,m-1),
1\ls i_m\ls i_{m-1}-1\},
\endgather
$$
and put
$$
A(n,m)=\sum_{(i_1,\dots,i_m)\in\Cal I(n,m)}g_n(i_1,\dots,i_m).
$$

\proclaim{Theorem 2.3} Let $n\gs d+2$ and $2\ls m\ls \min\{n, k-d\}-1$, then
$$
A(n,m)=A(n-1,m-1)+A(n,m+1).
$$
\endproclaim

\demo{Proof} By the definition,
$$
A(n,m)=\sum_{(i_1,\dots,i_m)\in\Cal I(n,m)}\sum_{j=1}^ng_n(i_1,\dots,i_m,j).
$$
Let us decompose the internal sum into three subsums over the intervals
$1\ls j\ls i_m-1$, $i_m\ls j\ls n-1$, $j=n$, and denote them $S_1$, $S_2$,
and $S_3$, respectively.
Evidently,
$$\align
S_1
&=\sum_{(i_1,\dots,i_m)\in\Cal I(n,m)}\sum_{j=1}^{i_m-1}g_n(i_1,\dots,i_m,j)\\
&=\sum_{(i_1,\dots,i_m,j)\in\Cal K(n,m+1)}g_n(i_1,\dots,i_m,j).
\endalign
$$
It is easy to see that if $(i_1,\dots,i_m,j)\in\Cal K(n,m+1)$ and 
$(i_1,\dots,i_m,j)\notin\Cal I(n,m+1)$ then there exists $p$, $1\ls p\ls m$,
such that $i_p=j$, and therefore by Lemma~2.2(i),  
$$
S_1=\sum_{(i_1,\dots,i_m,j)\in\Cal I(n,m+1)}g_n(i_1,\dots,i_m,j)=A(n,m+1).
$$

Next, 
$$
S_2=\sum_{(i_1,\dots,i_m)\in\Cal I(n,m)}\sum_{j=i_m}^{n-1}g_n(i_1,\dots,i_m,j)=
0
$$ 
by Lemma~2.2(i) and~(ii). Finally, by Lemma~2.2(iii),
$$\align
S_3&=\sum_{(i_1,\dots,i_m)\in\Cal I(n,m)}g_n(i_1,\dots,i_m,n)\\
&=\sum_{(i_1,\dots,i_m)\in\Cal J(n-1,m)}g_{n-1}(i_1,\dots,i_m).
\endalign
$$
Denote the latter expression by $B(n-1,m)$, so 
$$
A(n,m)=A(n,m+1)+B(n-1,m),  \tag2
$$
for $n\gs d+1$, $1\ls m\ls \min\{n,k-d\}-1$.

 Consider now $B(n,m+1)$. One has
$$\align
B(n,m+1)&=\sum_{(i_1,\dots,i_m,j)\in\Cal I(n,m+1)}g_n(i_1,\dots,i_m,j)\\
&+\sum_{(i_2-1,\dots,i_m-1,j-1)\in\Cal J(n,m)}g_n(n-d+1,i_2,\dots,i_m,j),
\endalign
$$
and by Lemma~2.2(iv), the second term in the right hand side equals $B(n-1,m)$.
Therefore,
$$
B(n,m+1)=A(n,m+1)+B(n-1,m) 
$$
for $n\gs d+1$, $1\ls m\ls \min\{n,k-d\}-1$.

Comparing this with (2) one gets $B(n,m+1)=A(n,m)$
for $ n\gs d+1$, $1\ls m\ls \min\{n,k-d\}-1$, and hence
$$
A(n,m)=A(n-1,m-1)+A(n,m+1)
$$
for $n\gs d+2$, $2\ls m\ls \min\{n,k-d\}-1$.
\qed
\enddemo

Our next aim is to express $A(n,m)$ in terms of $f_j$'s. We start from the 
following result.

\proclaim{Lemma~2.4} Let $n\gs d+2$ and $1\ls j\ls d$, then
$$
g_n(n+1-j)=\sum_{i=0}^{j-1}(-1)^i\binom{j-1-i}i f_{n-1-i}. \tag3
$$
\endproclaim

\demo{Proof} Let $n\gs3$; it is clear that $g_n(n)=f_{n-1}$, hence (3) 
holds for $j=1$ and $n\gs 3$.
Moreover, if $d\ne1$ then $g_n(n-1)=f_{n-1}$ as well, and hence (3) holds
also for $j=2$ and $n\gs 3$. This gives the basis for the induction.

Next, observe that
$$
g_n(n+2-j)-g_n(n+1-j)=g_{n-1}(n+2-j)
$$
for $n\gs 3$, $3\ls j\ls \min\{d,n\}$. Indeed, 
$$
g_n(r)=\sum_{i=1}^ng_n(r,i).
$$
The sum in the right hand side is subdivided into four sums: for 
$1\ls i\ls r-1$, for $i=r$, for $r+1\ls i \ls n-1$, and for $i=n$,
respectively. Using
Lemma~2.2(iv), (i), (ii), and (iii) respectively, we get
$$
g_n(r)=\sum_{i=1}^{r}g_{n-1}(i)
$$
for $n\gs 3$ and $n+1-d\ls r\ls n-1$. Observe that condition 
$m\ls  k-d-1$ from Lemma~2.2(iii) is satisfied by virtue of (1)
and $k\gs3$.
It remains to substitute $r=n-j+1$ and $r=n-j+2$, 
and to find the difference between the right hand sides; the restrictions 
on $r$ are transformed into $3\ls j\ls d$.

Finally, it is easy to check that
$$\multline
\sum_{i=0}^{j-2}(-1)^i\binom{j-2-i}i f_{n-1-i}-
\sum_{i=0}^{j-3}(-1)^i\binom{j-3-i}i f_{n-2-i}\\=
\sum_{i=0}^{j-1}(-1)^i\binom{j-1-i}i f_{n-1-i},
\endmultline
$$
for any $n\gs2$ and $j\gs3$, 
and hence (3) is satisfied for $n\gs 3$ and $1\ls j\ls \min\{d, n-1\}$.
It follows now from (1) that (3) is satisfied for $n\gs d+2$ and 
$1\ls j\ls d$.
\qed
\enddemo

We are now ready to find an expression for $A(n,m)$ similar to expression
(3) for $g_n(n+1-j)$.

\proclaim{Theorem 2.5} Let $n\gs k$ and $1\ls m\ls k-d$, then
$$
A(n,m)=\sum_{i=0}^{m+d}(-1)^i\binom{m+d-i}i f_{n-i}. \tag4
$$
\endproclaim

\demo{Proof} Similarly to the proof of Lemma~2.4, we first check (4)
for $m=1$ and $m=2$. Clearly,
$$
A(n,1)=\sum_{j=1}^{n-d}g_n(j)=f_n-\sum_{j=1}^dg_n(n+1-j).
$$
Using (3) we get
$$\align
A(n,1)&=f_n-\sum_{j=1}^d\sum_{i=0}^{j-1}(-1)^i\binom{j-1-i}i f_{n-1-i}\\
      &=f_n-\sum_{i=0}^{d-1}(-1)^i f_{n-1-i}\sum_{j=i+1}^d\binom{j-1-i}i\\
      &=\sum_{i=0}^{d+1}(-1)^i\binom{d+1-i}i f_{n-i}
\endalign
$$
for $n\gs d+2$, so (4) holds for $m=1$ and $n\gs d+2$.

Next, by (2), we have
$$
A(n,2)=A(n,1)-B(n-1,1)=A(n,1)-A(n-1,1)-g_{n-1}(n-d)
$$
for $n\gs d+2$ (condition $m=1\ls k-d-1$ follows from (1) and $k\gs3$ as
in the proof of Lemma~2.4). Using (3) for $j=d$ and (4) for $m=1$ we get
$$\multline
A(n,2)= \sum_{i=0}^{d+1}(-1)^i\binom{d+1-i}i f_{n-i}\\-
\sum_{i=0}^{d+1}(-1)^i\binom{d+1-i}i f_{n-1-i}-
\sum_{i=0}^{d-1}(-1)^i\binom{d-1-i}i f_{n-2-i}\\
=\sum_{i=0}^{d+2}(-1)^i\binom{d+2-i}i f_{n-i}
\endmultline
$$
for $n\gs d+2$, so (4) holds for $m=2$ and $n\gs d+2$.

Finally, by Theorem~2.3,
$$
A(n,m)=A(n,m-1)-A(n-1,m-2)
$$
for $n\gs d+2$ and $3\ls m\ls k-d$, so by induction we get
$$\align
A(n,m)&=\sum_{i=0}^{m-1+d}(-1)^i\binom{m-1+d-i}i f_{n-i}-
\sum_{i=0}^{m-2+d}(-1)^i\binom{m-2+d-i}i f_{n-1-i}\\
&=\sum_{i=0}^{m+d}(-1)^i\binom{m+d-i}i f_{n-i}
\endalign
$$
for $n\gs d+2$ and $1\ls m\ls\min\{n+1,k\}-d$, and hence for
$n\gs k$ and $1\ls m\ls k-d$.
\qed
\enddemo

Now we are ready to prove Theorem~1.2 for $k\gs3$. First of all,
it follows from Lemma~2.2(v) that 
$A(n,k-d)=0$ for all $n\gs k$. Hence, by Theorem~2.5,
$$
\sum_{i=0}^k(-1)^i\binom{k-i}if_{n-i}=0
$$
for $n\gs k$, or, equivalently,
$$
\sum_{i=0}^k\sum_{n\gs k}(-1)^ix^i\binom{k-i}if_{n-i}x^{n-i}=0.
$$

Evidently, for $n\ls k-1$ the layered restriction is void, so by 
\cite{Kn, SS},
$f_n=c_n$ for $0\ls n\ls k-1$, where $c_n$ is the $n$th Catalan number. 
Therefore
$$
\sum_{i=0}^k(-x)^i\binom{k-i}i\left(F_T(x)-\sum_{j=0}^{k-i-1}x^jc_j\right)
=0. \tag5
$$

Recall that the Chebyshev polynomials of the second kind satisfy relation
$$
x^{k/2}U_k\left(\frac 1{2\sqrt{x}}\right)=\sum_{i=0}^k(-x)^i\binom{k-i}i
$$
(see \cite{Ri, pp.~75,76}), while the Catalan numbers satisfy relation
$$
\sum_{i=0}^l(-1)^i\binom{k-i}ic_{l-i}=(-1)^l\binom{k-1-l}l
$$
for $l\ls k-1$
(see \cite{Ri, pp.~152--154}), and hence
$$
\sum_{i=0}^k(-x)^i\binom{k-i}i \sum_{j=0}^{k-i-1}x^jc_j=
\sum_{l=0}^{k-1}x^l\sum_{i=0}^l(-1)^i\binom{k-i}ic_{l-i}=
\sum_{l=0}^{k-1}(-x)^l\binom{k-1-l}l.
$$
Therefore, (5) yields
$$
F_T(x)=\frac{\sum_{i=0}^{k-1}(-x)^i\binom{k-1-i}i}
{\sum_{i=0}^{k}(-x)^i\binom{k-i}i}=
\frac{U_{k-1}\left(\frac 1{2\sqrt{x}}\right)} 
{\sqrt{x}U_{k}\left(\frac 1{2\sqrt{x}}\right)}=R_k(x)
$$
for $k\gs 3$.

The remaining case $k=2$ is trivial, since $R_2(x)=(1-x)^{-1}$, which is 
equivalent to the fact that there is exactly one permutation of length $n$
avoiding $w(1,1)=(1,2)$.

\enddocument